\documentclass[12pt]{article}
\input epsf
\usepackage{amssymb}

\def \be{\begin{equation}}
\def \ee{\end{equation}}
\def \berr{\begin{eqnarray}}
\def \err{\end{eqnarray}}
\def \nn{\nonumber}
\def \a{\alpha}
\def \ad{\alpha^{\vee}}
\def \b{\beta}

\def \d{\delta}

\def \Del{\Delta}

\def \l{\lambda}
\def \L{\Lambda}

\def\c#1{{\cal #1}}

\def \({\left(}
\def \){\right)}
\def \<{\langle}
\def \>{\rangle}
\def \[{\left[}
\def \]{\right]}

\def\tens{\mathop{\otimes}}

\def\reps{representations }
\def\rep{representation }

\def\mf#1{{\mathbb #1}}
\def\Z{{\mf{Z}}}
\def\N{{\mf{N}}}
\def\Q{{\mf{Q}}}
\def\R{{\mf{R}}}
\def\C{{\mf{C}}}

\def\mg{\mathfrak{g}}

\newcommand{\sect}[1]{\setcounter{equation}{0}\section{#1}}

\newtheorem{prop}{Proposition}[section]
\newtheorem{theorem}[prop]{Theorem}

\evensidemargin \oddsidemargin
\begin{document}

\begin{titlepage}

\begin{center}  

{\Large\bf
Unitary Representations of Noncompact Quantum Groups 
at Roots of Unity
 \\[4ex]}
H. \ Steinacker 
\footnote{Harold.Steinacker@physik.uni--muenchen.de}
\\[2ex] 
{\small\it 
          Institut f\"ur Theoretische Physik \\        
          Ludwig--Maximilians--Universit\"at M\"unchen\\
        Theresienstr.\ 37, D-80333 M\"unchen  \\[1ex] }
\end{center}

\vspace{5ex}

\begin{abstract}

Noncompact forms of the Drinfeld--Jimbo quantum groups $U_q^{fin}(\mg)$ 
with $H_i^* = H_i$, \\ ${X_i^{\pm}}^* = s_i X_i^{\mp}$ for $s_i=\pm 1$ 
are studied at roots of unity. This covers $\mg = so(n,2p),su(n,p)$, 
$so^*(2l)$, $sp(n,p)$, $sp(l,\R)$, and exceptional cases. 
Finite--dimensional unitary \reps are found for all these forms,
for even roots of unity. Their classical symmetry induced by the 
Frobenius--map is determined, and the meaning of the extra quasi--classical 
generators appearing at even roots of unity is clarified. 
The unitary highest weight modules of the classical case 
are recovered in the limit $q \rightarrow 1$.
\end{abstract}

\parskip 4pt plus2pt minus2pt
\vfill
\noindent
Preprint LMU-TPW 99-09

\end{titlepage}

\sect{Introduction}

Quantum groups allow to generalize the concept of symmetry, 
which has proved to be of great importance in physics. 
Up to this date, most of the work on quantum groups
has been done for the 
compact case. However noncompact groups are important as well,
for example the Lorentz group, or the Anti--de Sitter group $SO(2,n)$  
which has attracted much attention recently in the 
context of string theory \cite{maldacena}. 

We consider the Drinfeld--Jimbo quantized universal enveloping algebra
$U_q^{res}(\mg)$ \cite{drinfeld,FRT,jimbo}
corresponding to finite--dimensional semisimple Lie algebras.
In the $q$--deformed case, there are several possibilities to define 
real, in particular noncompact forms of these algebras. 
If $q$ is real, the \rep theory is largely parallel to the 
classical case, but more complicated; for some results  
in this case see \cite{klimyk,guizzi}.
In the present paper, we consider instead the case where $q$ 
is a root of unity, which provides additional structure 
that does not exist in the classical case. This
turns out to be much simpler, rather than more difficult than
the undeformed case. We study unitary \reps of 
(a slight extension of) the so--called ``finite'' quantum group 
$U_q^{fin}(\mg) \subset U_q^{res}(\mg)$ at roots of unity, 
with real structure of the form $H_i^* = H_i$ and 
$X_i^{\pm *} = s_i X_i^{\mp}$, where $s_i =\pm 1$. 
This covers  $so(n,2p),su(n,p), so^*(2l), sp(n,p)$, $sp(l,\R)$, 
as well as various forms for the exceptional groups. 
Even though this real form corresponds to a non--standard
Hopf algebra $*$--structure, it is appropriate for our 
purpose, and leads to a large class of unitary representations.

Generalizing the method of \cite{ads_paper}, 
we find unitary \reps for all these noncompact forms,
provided $q$ is an even roots of unity. It is shown that
all of them can be related to unitary \reps 
of the compact form in a simple way. As opposed to the classical case,
they are finite--dimensional, 
which means that the problem is a purely algebraic one. In many cases,
they can be viewed as regularizations of classical, 
infinite--dimensional representations. In particular, we show how 
almost all classical unitary highest weight modules 
(with the possible exception of a certain ``small'', 
discrete set of highest weights)
can be obtained as the limit $q\rightarrow 1$
of unitary \reps of $U_q^{fin}$.
In the example of the Anti--de Sitter group $SO(2,3)$, this was 
already studied for special cases in \cite{dobrev,flato}, and
more generally in \cite{ads_paper}. 
Not all the \reps found however 
have a classical limit in an obvious way; to understand this better
is an interesting open problem. Moreover, 
it turns out that the unitary \reps of 
$U_q^{fin} \subset U_q^{res}(\mg)$ are very different 
from the ones studied in \cite{schmuedgen}, where a different 
specialization of $U_q(sl(2,\R))$ to roots of unity is considered,
leading to an infinite--dimensional algebra.

This paper is organized as follows. After reviewing the definitions
and basic concepts in Section 2, the unitary \reps of the compact case are 
studied in Section 3, and the particular features appearing at roots 
of unity are discussed. In Section 4, the remarkable classical symmetry
$U(\tilde \mg)$ arising from $U_q^{res}(\mg)$ at roots of unity due to the 
Frobenius map \cite{lusztig_90,lusztig_book} is discussed, including the
case of even roots of unity which turns out to be most important. 
The extra generators arising at even roots of unity
which extend the classical universal enveloping algebra 
find a natural interpretation here. 

In Section 5, the noncompact forms are defined, and unitary \reps
are found for all of them in a rather simple way. It turns out that 
only a subgroup of the classical $U(\tilde \mg)$ preserves the 
noncompact form, which is determined in Section 6. Finally in Section 7,
the connection with the classical case is made, and it is shown how 
the classical unitary highest--weight \reps are recovered in the limit 
$q \rightarrow 1$.

In the appendix, 
an explicit, self--contained approach to the classical symmetry
arising from the Frobenius map is given 
including the case of even roots of unity, 
which was treated only implicitly in \cite{lusztig_book}.

\sect{Definitions and basic properties}

We first collect the basic definitions, in order to fix the notation.
Let $A_{ij}=2\frac{(\a_i, \a_j)}{(\a_j,\a_j)}$ be
the Cartan matrix of a classical simple Lie algebra $\mg$ of rank $r$,
where $( , )$ is the inner product in root space 
and $\{\a_i,\;  i=1,...,r \} $
are the simple roots. The positive roots will be denoted by $Q^+$,
and $\rho = \frac 12 \sum_{\a \in Q^+} \a$ is the Weyl vector.

For $q\in \C$, the {\em quantized universal enveloping algebra}
$U_q(\mg)$ is the Hopf algebra with generators
$\{ Y^{\pm}_i, K_i, K_i^{-1}; \; i=1,...,r\}$ and 
relations \cite{jimbo,drinfeld,FRT}
\berr
\[K_i, K_j\]      &=& 0 ,                   \\  
 K_i Y^{\pm}_j    &=& q^{\pm A_{ji}} Y^{\pm}_j K_i,     \\
\[Y^+_i, Y^-_j\]   &=& \d_{i,j} \frac{K_i
                          -K_i^{-1}}{q^{d_i}-q^{-d_i}}, \\
\sum_{k=0}^{1-A_{ji}} &\ &
           \[\begin{array}{c} 1-A_{ji} \\ k\end{array}\]_{q_i}(Y^{\pm}_i)^k
         Y^{\pm}_j (Y^{\pm}_i)^{1-A_{ji}-k} = 0, \quad i\neq j,
\label{UEA_K}
\err
where the $d_i = (\a_i, \a_i)/2$ are relatively prime, 
$q_i = q^{d_i}, \quad [n]_{q_i} = \frac{q_i^n-q_i^{-n}}{q_i-q_i^{-1}}$, and
\be
\[ \begin{array}{c} n\\m \end{array} \]_{q_i} =
             \frac{[n]_{q_i}!}{[m]_{q_i}! [n-m]_{q_i} !} .  \label{binom}
\ee
We assume that $q^{d_i}\neq q^{-d_i}$.
The comultiplication is defined by 
\berr
\Del(K_i)          &=& K_i \tens K_i,  \nonumber \\
\Del(Y^{+}_i) &=&  1\tens Y^{+}_i  +  Y^{+}_i\tens K_i,   \nonumber \\
\Del(Y^{-}_i) &=&  K_i^{-1} \tens Y^{-}_i  + Y^{-}_i \tens 1.
\label{coproduct_Y}
\err
Antipode and counit exist as well, but will not be needed. 
The Borel subalgebras $U_q^{\pm}(\mg)$
are defined in the obvious way.

In this paper, $q$ will always be a complex number, rather than a formal
variable. Moreover, since we are mainly interested in 
representations, it is more intuitive
to use the generators $\{X_i^{\pm}, H_i\}$ defined by
\be
K_i = q^{d_i H_i},  \quad   Y_i^{+} = X_i^{+} q^{H_i d_i/2}, 
\quad Y_i^{-} = q^{-H_i d_i/2} X_i^-,
\label{XY_def}
\ee
so that the relations take the more familiar form
\berr
\[H_i, H_j\]      &=& 0,                      \label{CR_HX} \\
\[H_i, X^{\pm}_j\] &=& \pm A_{ji} X^{\pm}_j,     \\
\[X^+_i, X^-_j\]   &=& \d_{i,j} \frac{q^{d_i H_i}
                          -q^{-d_i H_i}}{q^{d_i}-q^{-d_i}}
                    = \d_{i,j} [H_i]_{q_i}.  
\label{UEA}
\err
The comultiplication is now
\berr
\Del(H_i)          &=& H_i \tens 1 + 1 \tens H_i,       \nonumber \\
\Del(X^{\pm}_i) &=&  X^{\pm}_i \tens q^{d_i H_i/2} + q^{-d_i H_i/2}
\tens  X^{\pm}_i.
\label{coproduct_X}
\err
The classical case is recovered for $q=1$.
Generators $X_{\a}^{\pm}$ 
corresponding to the other positive roots $\a$ can be defined using
the braid group action \cite{lusztig}; we will quote some 
properties as they are needed.
A Poincare--Birkhoff--Witt (P.B.W.) basis 
is then given as classically in terms of ordered  monomials 
of the raising and lowering operators corresponding to all positive 
respectively negative roots.

If $q$ is allowed to be a root of unity, we will instead consider the 
``restricted specialization'' $U_q^{res}(\mg)$ \cite{lusztig}
with generators $X_i^{\pm(k)} = \frac{(X_i^{\pm})^k}{[k]_{q_i}!}$ 
for $k \in \N$ as well as $H_i$. For generic $q$,
i.e. $q$ not a root of unity, this is the same as before.
However if $q$ is a root of unity, 
\be
q=e^{2\pi i n/m}   \label{q_root}
\ee
with $m$ and $n$ relatively prime,  then $[k]_q$ becomes 0 for certain $k$.
Denote with $M$
the smallest positive integer such that $q^{2M} =1$, i.e. $[M]_q =0$. Thus
$M=m$ if $m$ is odd, and $M=m/2$ if $m$ is even. 
In the first case $q^M=1$, and we will say that $q$ is an ``odd'' 
root of unity. In  the second case $q^M =-1$, and $q$ will be called ``even''.
More generally for $q_i= e^{2\pi i d_i n/m}$, let $M_i$ be the  
the smallest integer such that
\be
[M_i]_{q_i} =0.
\ee
Then $M_i$ divides $M$; similarly, we define $M_{\a}$ and $d_{\a}$ 
for the other roots. 
$U_q^{res}(\mg)$ contains the additional generators $X_i^{\pm (M_i)}$,
which have a well--defined coproduct, and thus are defined on
tensor products of representations. Verma modules can also be defined  
in the usual way, for integral highest weights \cite{lusztig_mod}.
We will only consider these types of representations of 
$U_q^{res}(\mg)$. In particular,
$(X_i^{\pm})^{M_i}=0$ in $U_q^{res}(\mg)$. Therefore
$U_q^{res}(\mg)$ contains a remarkable sub--Hopf algebra
$u_q^{fin}$ (the ``small quantum group'') generated by
$X_i^{\pm}$ and  $K_i^{\pm 1}$.
We prefer to slightly change the standard convention and 
define $U_q^{fin}$ by including the $H_i$ as well,
slightly abusing the name ``finite''. 
This is a more intuitive generalization of the classical 
$U(\mg)$ at least from a physical point of view, and poses no
problem since $q$ is a complex number here rather than a 
formal variable.

The generators $X_i^{\pm (M_i)}$  act 
as (graded) derivations on $u_q^{fin}$ by 
$x \rightarrow [X_i^{\pm(M_i)}, x]_{\pm}$.
The right--hand side is indeed an element of $u_q^{fin}$, as can be seen 
from the commutation relations (\ref{XX_cr}). 

Finally, we quote the following useful relation:
\be
[X_i^+, (X_i^-)^k] = (X_i^-)^{k-1} [k]_{q_i} [H_i-k+1]_{q_i}.
\label{keller_eq}
\ee

\sect{Representations of $U_q^{fin}$ and weight space}

The Cartan generators can be evaluated on weights $\l$, such that
\be
<H_i,\l> = \frac {(\a_i,\l)}{d_i} = (\ad_i,\l),
\ee
where as usual $\ad = \frac{2\a}{(\a,\a)}$ 
is the coroot of $\a$.
The fundamental weights $\L_i$ satisfy $(\L_i, \ad_j) =  \d_{i,j}$,
therefore
\be
<H_i,\L_j> = \d_{ij},
\ee
and span the lattice of integral weights. The Weyl group $\c{W}$ 
is defined as usual, and 
$D = \{\sum_i r_i \L_i; \; r_i \in \R_{\geq 0}\}$ is the dominant Weyl
chamber.

It is well--known \cite{rosso} that for generic 
$q$, the representation theory is 
essentially the same as in the classical case. In particular, the 
finite--dimensional \reps (=modules)  of $U^{res}_q(\mg)$ are direct sums of
irreducible representations  $L^{res}(\l)$, which are highest--weight
\reps with dominant integral highest weight $\l$. Their character 
\be
\chi (L^{res}(\l)) = e^{\l}\sum_{\eta >0}\dim{L^{res}(\l)_{\eta}}e^{-\eta}
  =:\chi(\l)
\label{char}
\ee
is given by Weyls formula. Here $L^{res}(\l)_{\eta}$ is the weight space
of $L^{res}(\l)$ with weight $\l-\eta$.
Irreducible highest weight \reps of $U_q^{fin}$ are denoted by $L^{fin}(\l)$.

\subsection{Singlets, special points, and the dual algebra $\tilde \mg$}

One important feature at roots of unity is the existence of 
nontrivial one--dimensional 
\reps $L^{fin}(\l_z)$ of $U_q^{fin}$, with weights 
\be
\l_{z} = \sum_i z_i M_i \L_{i}
\label{l_z}
\ee
for $z_i \in \Z$; this follows from  (\ref{keller_eq}).
There also exist similar \reps with $z_i \notin \Z$ which will be
considered in Section 5, but for now  we concentrate on the case 
of integral weights. These weights $\l_{z}$
will be called {\em special points}. They span a lattice which is 
the weight lattice of a dual Lie algebra $\tilde \mg$, 
rescaled by $M$. 
In particular, it contains the root lattice of $\tilde \mg$, which is
generated by the $M_i \a_i$ or equivalently $M_{\a} \a$. 
Indeed, consider a second metric on 
weight space defined by \cite{lusztig_book}
\be
(\a_i,\a_j)_d:=(M_i \a_i,M_j \a_j),             \label{dual_metric}
\ee
with associated matrix
\be
\tilde A_{ij}:= 2\frac{(\a_i,\a_j)_d}{(\a_j,\a_j)_d} = \frac{M_i}{M_j} A_{ij}.
\ee
In particular,
\be
M_i A_{ij} = \tilde A_{ij} M_j.                    \label{dual_cartan}
\ee
$\tilde A_{ij}$ is always a Cartan Matrix: it is clearly nondegenerate,
and $\tilde A_{ii} =2$. To see that $\tilde A_{ij} \in -\N_0$ 
for $i\neq j$, observe that by the definition of $M_j$,
$M_j d_j$ is the smallest integer which is 
divisible by both $M$ and $d_j$.
Similarly $A_{ji} d_i M_i$ is divisible by $M$ because $A_{ji}$ is an 
integer, and also by $d_j$, since $A_{ji} d_i = A_{ij} d_j$. Therefore 
$\frac{A_{ji} d_i M_i}{M_j d_j}$ is an integer, equal to
$\frac{M_i}{M_j} A_{ij} = \tilde A_{ij}$.

We shall determine $\tilde A_{ij}$ explicitely. 
In the simply laced  case,
all $M_i$ are equal, therefore $(\;,\;)_d$ is proportional to the
Killing metric, and $\tilde A_{ij} =  A_{ij}$. Thus the lattice
of special points is nothing but the weight lattice rescaled by $M$,
and $\tilde \mg = \mg$.

For $B_n, C_n$ and  $F_4$, there are roots with 2 different lengths 
$d_s=1$ and $d_l=2$.
Again, if $M$ is not divisible by 2, i.e. if $q$ is odd, then clearly 
$M_i = M$ is odd for all $i$, and $\tilde A_{ij} =  A_{ij}$.
On the other hand  if $q$ is even, 
then $M_i = M/2 =:M_l$ if $\a_i$ is long, and
$M_i = M =: M_s$  if $\a_i$ is short.
Thus $\frac{M_i}{M_j} = \frac{d_j}{d_i}$ and $\tilde A_{ij} = A_{ji}$,
which means that $M_l \a_l$ are the short roots and
$M_s \a_s$ the long roots in the lattice of special points.
Therefore the dual algebra of $B_n$ is $C_n$ and vice versa, 
while $F_4$ remains $F_4$ except that the roots change their role.

For $G_2$, the roots have lengths $d_s=1$ and $d_l=3$.
If $M$ is not divisible by 3,
then  $M_i = M$ for $i=1,2$, and again $\tilde A_{ij} =  A_{ij}$.
On the other hand if $M$ is divisible by 3, let $\a_l:=\a_1$ be the
long simple root, and $\a_s:=\a_2$ be the short one. Then 
$M_l = M_1 = M/3$, $M_s = M_2 = M$, and $\tilde A_{ij} =  A_{ji}$.
Thus the dual lattice is again of type $G_2$, but now $M_s \a_s$ is the
long root, and $M_l \a_l$ the short one.

To summarize, $\tilde \mg = \mg$, except for $\tilde{B_n} = C_n$
and $\tilde{C_n} = B_n$ if $q$ is even.
For all cases, the Weyl group of $\tilde \mg$ is the same 
as that of $\mg$. In Section 4, we will see that in some sense,
$U_q^{res}(\mg)$ contains indeed a classical algebra associated 
with the lattice of special points.

The hyperplanes 
\be
H_{\a}^z:=\{\l; \; (\l, \ad) = M_{\a} z \},
\label{H_z}
\ee
where $\a$ is any root and $z \in \Z$, 
divide weight space into simplices called {\em alcoves}. 
The alcove of dominant weights with the origin on its boundary
is called the {\em fundamental alcove}. 
The reflections on these hyperplanes generate the 
affine Weyl group, which plays an important role in the 
\rep theory at roots of unity. 
Notice that every special point is in some $H_{\a}^z$ for every root $\a$.
To see this, we have to show that $(M_i\L_i, \ad) \in M_{\a} \Z$
for every root $\a$. Since the Weyl group preserves the lattice
generated by $M_i\L_i$, this follows from  the fact that
$(\sum_i z_i M_i\L_i,\ad_j) \in M_j \Z$,
for a suitable $\a_j$. In fact, the special points
are the intersection points of a maximal number of hyperplanes.

\subsection{Unitary representations of the compact form}

To define unitary representations, one first has to 
specify the real form of the algebra, or group in the classical case. 

A real form or $\ast$ -- structure  
is an antilinear involution (=anti--algebra map) on $U_q^{res}(\mg)$.
In the classical case, the $\ast$ is acting on the complexified 
Lie algebra, and the real Lie algebra is by definition its
eigenspace with eigenvalue $-1$. 
The interpretation of a real form at
$q\neq 1$ is given by its classical limit. 

In this section, we only consider the compact form. It is defined
by $\ast = \theta$ where 
$\theta(X_i^{\pm}) = X_i^{\mp},\; \theta(H_i) = H_i$ 
is the Cartan--Weyl involution, thus
\be
(X_i^{\pm})^* = X_i^{\mp}, \quad  H_i^* = H_i,
\label{formal_invol}
\ee
extended as an antilinear anti--algebra map. This is consistent 
for $q$ real and $|q|=1$.

A representation of $U_q^{res}(\mg)$ on a Hilbert space $V$ is said to be
{\em unitary} if 
the star is implemented as the adjoint on the Hilbert space,
i.e. $(v,x\cdot w) = (x^{\ast}\cdot v,w)$ for any 
$x \in U_q^{res}(\mg)$ and $v,w \in V$. In particular, $(\;,\;)$
is positive definite. In the classical case, this means precisely 
that the adjoint (=star) of a group element is its inverse.
Since all unitary \reps are completely reducible,
we only need to consider irreducible ones.
Unitary and unitarizable will be used synonymously.
On unitary highest weight modules with (\ref{formal_invol}) 
or (\ref{s_i_reality}), the inner product 
can be calculated recursively, descending from the highest weight state.
In particular, it is unique up to normalization. 

Finite--dimensional unitary \reps
of noncompact forms with the correct classical limit
are possible only at roots of unity. 
Therefore we will concentrate on that case from now on, 
in particular $q^{\ast} = q^{-1}$. 
Even though (\ref{formal_invol}) is then a ``nonstandard''
Hopf algebra $\ast$--structure\footnote{since $(S(x))^* = S(x^*)$ here
where $S$ is the antipode,
rather than $(S(x))^* = S^{-1}(x^*)$ for $x \in U_q^{res}(\mg)$}, 
it is appropriate for our purpose.

All finite--dimensional \reps of $U_q^{res}(\mg)$ 
have integral weights, even at roots of unity.
While this is not true for $U_q^{fin}$
any more, we nevertheless start with studying the unitary
\reps of $U_q^{fin}$ with integral weights. 
The following well--known fact \cite{CH_P} is useful:
\begin{theorem}
Assume that  $\l$ is a dominant integral weight with
$(\l + \rho, \ad) \leq M_{\a}$ for all positive roots $\a$.
Then the highest weight \rep $L^{res}(\l)$ has the same character 
$\chi$ as in the 
classical case, given by Weyl's character formula.
\label{char_thm}
\end{theorem}
In other words, $\l+\rho$ is in the fundamental alcove. 
This follows from the strong linkage principle, 
which was first shown in \cite{anderson}; 
for a more elementary approach, see \cite{ads_paper}.
Moreover, $L^{fin}(\l) = L^{res}(\l)$
for these weights $\l$, since the $X_{i}^{\pm (M_{i})}$ act trivially.
If the above bound is not satisfied, then the Verma module with highest 
weight $\l$ contains additional highest weight submodules 
besides the classical ones.

Now we can show the following:
\begin{theorem}
Let $\l$ be a dominant integral weight, and  $q=e^{2\pi i n/m}$.
Then $L^{fin}(\l)$ is a unitary \rep of the compact form 
(\ref{formal_invol}) of $U_q^{fin}$
if the character of $L^{fin}(\l)$ is given by Weyl's formula
for all $q'=e^{2\pi i \varphi'}$ with $0\leq \varphi'<n/m$.
In particular, this holds if 
$(\l + \rho,\ad) \leq \lceil\frac{m}{2nd_{\a}}\rceil+ 1$ 
for all positive roots $\a$, where $\lceil c \rceil$ denotes the 
largest integer $\leq c$ for $c\in \R$.
\label{cpct_thm}
\end{theorem}
\begin{proof}
Consider $L^{fin}(\l)$
for all $q' \in B:=\{e^{2\pi i \varphi};\;  0 \leq \varphi < n/m\}$.
If the character of $L^{fin}(\l)$ is the same for all $q\in B$,
one can identify the
$L^{fin}(\l)$ as vector spaces\footnote{or even better, 
view them as trivial vector bundle
over $B$, with local trivializations given in terms of the P.B.W. basis}. 
Their inner product matrix is smooth (in fact analytic) 
in $q'$, and positive definite at $q'=1$ since we consider the
compact case. This implies that all eigenvalues are positive on $B$:
assume to the contrary that the matrix were not positive definite 
for some $q' \in B$. Then it would have a zero eigenvalue for some
$q_0 \in B$, which implies that its null space is a 
submodule of $L^{fin}(\l)$. But this is impossible, 
since the $L^{fin}(\l)$ are irreducible by definition.
For $q'=e^{2\pi i n/m}$, some eigenvalues may vanish; but then
$L^{fin}(\l)$ is the quotient of 
$\lim_{q'\rightarrow q \atop q' \in B} L^{fin}(\l)$ 
modulo its null space, which again has a 
positive definite inner product.

In particular, assume that 
$(\l + \rho, \ad) \leq \lceil\frac{m}{2nd_{\a}}\rceil+ 1$.
Let $M'_{\a}$ be the smallest integer $> \frac m{2d_{\a}n}$,
which is $\lceil\frac{m}{2nd_{\a}}\rceil+ 1$.
Then $M'_{\a}$ is associated to 
$q':=e^{\frac{2\pi i}{2d_{\a}M'_{\a}}} \in B$ as defined 
in Section 2. Therefore by  Theorem \ref{char_thm}, 
$L^{fin}(\l)$ at $q'$  has the same character as for $q=1$, since 
$(\l + \rho, \ad) \leq \lceil\frac{m}{2nd_{\a}}\rceil+1 
= M'_{\a}$. For all other roots of unity $q'' \in B$, 
the character is again the same since the associated $M''_{\a}$ 
is larger than $M'_{\a}$. Thus the above argument applies.
\end{proof}

For some highest weights $\l$ on the boundary of the domain specified in 
Theorem \ref{cpct_thm}, the character of 
the unitary \rep $L^{fin}(\l)$ is smaller than the classical one.
The reason is that the generic \reps develop null--submodules;
this can be interpreted in the context of gauge theories, 
see \cite{ads_paper}. 

One may ask if all the unitary \reps have been found in Theorem 3.2.
As will be discussed in Section 7, it is possible
that there exist certain unitary \reps with integral weights
which do not even satisfy 
the first condition in Theorem \ref{cpct_thm}, 
as suggested by the classcial noncompact case. 
This would have to be studied by different methods. Other 
unitary \reps with integral and nonintegral weights 
will be obtained in Theorem \ref{nc_reps_thm}, 
which however do not have a classical limit.


\sect{Frobenius map and the quasi--classical symmetry $\tilde \mg$}

The modules $L^{res}(\l)=L^{fin}(\l)$ in Theorem \ref{cpct_thm} 
are irreducible \reps
of $U_q^{fin}$. For larger $\l$,  $L^{res}(\l)$ decomposes into a 
direct sum of irreducible modules of $U_q^{fin}$, 
which will be described now. 
This involves the special points introduced in Section 3.1.

The basic observation is the following.
Consider a highest--weight module $U_q^{-res}(\mg) \cdot v_{\l_z}$ 
with highest weight $\l_{z} = \sum_i z_i M_i \L_{i}$ and $z_i \in \Z$.
From  (\ref{keller_eq}), it follows that all $X_i^- \cdot v_{\l_z}$
are highest weight vectors. Therefore $X_i^- \cdot v_{\l_z}=0$ in 
$L^{res}(\l_z)$, because it is irreducible by definition. 
Using the P.B.W. basis, one can see that any element of 
$U_q^{-res}(\mg)$ can be written as a sum of terms of the form 
$(X_{\b_1}^{-(M_{\b_1})})^{k_1}...  (X_{\b_N}^{-(M_{\b_N})})^{k_N}
U_q^{-fin}$. It follows that all weights of 
$L^{res}(\l_z)$ have the form 
$\l_{z'} = \l_z - \sum_i n_i M_i \a_i$ with $n_i \in \N$.
In other words, $L^{res}(\l_{z})$ is a direct sum of 
one--dimensional \reps $L^{fin}(\l_{z'})$ of $U^{fin}_q$, since
$\l_{z'}$ is a special point.
However the ``large'' generators $X_{i}^{\pm (M_{i})}$ do act
nontrivially, as we will see.

Consider $L^{res}(\l_{z}) \tens L^{res}(\l_0)$ for $\l_z$ as above and
integral $\l_0$ with $0 \leq (\l_0,\ad_i) < M_i$.
Now the generators $Y^{\pm}_i$ (\ref{coproduct_Y}) are useful. 
Using the coproduct, one finds 
\be
Y_i^{\pm}\cdot(v\tens w) = v\tens (Y_i^{\pm} \cdot w),
\label{Y_1}
\ee
and 
\be
Y_i^{\pm (M_i)} \cdot (v\tens w) = (Y_i^{\pm (M_i)}\cdot v)
                                   \tens (K_i^{M_i}\cdot w)
\label{Y_2}
\ee
for $v \in L^{res}(\l_z)$ and $w \in L^{res}(\l_0)$,
because $Y_i^{\pm (M_i)} \cdot w =0$ by
the bound on $\l_0$. For the same reason, 
$L^{res}(\l_0)$ is an irreducible \rep of $U_q^{fin}$.
Together with (\ref{Y_1}) and (\ref{Y_2}), it follows that  
$L^{res}(\l_{z}) \tens L^{res}(\l_0)$ is an irreducible
highest weight module of $U_q^{res}(\mg)$, 
and we have verified \cite{lusztig_mod}

\begin{theorem}
Let  $\l_z$ and  $\l_0$ be integral weights as above with 
$0 \leq (\l_0,\ad_i) < M_i$, and $\l=\l_0+\l_z$. Then
\be
L^{res}(\l) = L^{res}(\l_0) \tens L^{res}(\l_{z}).
\label{L_res_l}
\ee
\label{L_res_theorem}
\end{theorem}
In particular, $L^{res}(\l)$ decomposes into  
a direct sum of irreducible \reps
$L^{fin}(\l - \sum_i n_i M_i \a_i)$ of  $U_q^{fin}$.
Moreover, (\ref{Y_1}) and (\ref{Y_2}) show that $Y_i^{\pm}$ commutes
with $Y_i^{\pm (M_i)} K_i^{M_i}$ on $L^{res}(\l)$. 

We will now see that the latter generators acting on $L^{res}(\l_z)$
provide a representation of the classical universal 
enveloping algebra $U(\tilde \mg)$
corresponding to the Cartan matrix $\tilde A_{ij}$.
This is the essence of
a remarkable result of Lusztig \cite{lusztig_90,lusztig_book}.
For odd roots of unity, it states 
that there is a surjective algebra homomorphism
\berr
U_q^{res}(\mg) &\rightarrow& U(\mg), \nn\\
Y_i^{\pm} &\rightarrow& 0  \nn\\
K_i    &\rightarrow& 1  \nn\\
Y_i^{\pm (M_i)} &\rightarrow& \tilde X_i^{\pm}.
\label{frobenius}
\err
This is the so--called  Frobenius map (recall that $\tilde \mg = \mg$
for odd roots of unity). 
It is generalized to even roots of unity 
in \cite{lusztig_book}; unfortunately the results given there are 
not very explicit. 
Since this case is of central importance to us, we will give an
elementary, self--contained approach, and show explicitely 
how the action of $U(\tilde \mg)$ on $L^{res}(\l_{z})$ is 
given in terms of the $X_j^{\pm (M_j)}$.
The complications arise because at even roots of unity,
$K_i$ cannot be set to 1, while $K_i^2$ must be, since 
$[Y_i^+,Y_i^-] = \frac {K_i-K_i^{-1}}{q_i-q_i^{-1}}$. Indeed, 
\be
<K_i,\sum_j z_j\L_j M_j> = q_i^{z_iM_i} = \pm 1.
\label{K_i_eval}
\ee
These extra, ``quasiclassical'' generators $K_i$
in some cases anticommute with $X_j^{\pm (M_j)}$, and 
will extend the algebra $U(\tilde \mg)$.
They will play an important role in the noncompact case.

Let $a_i \in \{0,1\}$ such that 
$a_i + a_j = 1$ if $\tilde A_{ij} \neq 0$ and $i\neq j$;
this is always possible.
Define $\tilde K_i = K_i^{M_i}$, and
\berr
\tilde X_i^+ &=& X_i^{+(M_i)} \tilde K_i^{a_i},\nn\\
\tilde X_i^- &=& X_i^{-(M_i)} \tilde K_i^{1-a_i}q_i^{M_i^2},\nn\\
\tilde H_i   &=& [\tilde X_i^+,\tilde X_i^-]
\label{class_generators}
\err
Then we can show the following:

\begin{theorem} 
For all special points $\l_{z}$, 
$L^{res}(\l_z)$ is an irreducible highest--weight \rep of the classical
$U(\tilde \mg)$, with generators
$\tilde X_i^{\pm}$ and $\tilde H_i$. 
If $v_{z'} \in L^{res}(\l_{z})$ has weight $\sum_j z'_j M_j \L_j$,
then $\tilde H_i \cdot v_{z'} = z'_i v_{z'}$.
Moreover, 
\be
\tilde X_i^{\pm} \tilde K_j = s_{ij} \tilde K_j \tilde X_i^{\pm} 
\ee
where 
$s_{ij} = q_i^{M_i M_j A_{ji}} = q^{M_i M_j (\a_i,\a_j)} = \pm 1$.
For dominant integral $\l$, 
$L^{res}(\l)$ is a direct sum of such irreducible representations,
by (\ref{Y_2}).
\label{classical_algebra}
\end{theorem}
This is proved in the appendix. Root vectors 
$\tilde X_{\tilde\a}^{\pm} \in U_q^{res}(\mg)$
for the remaining roots $\tilde\a \in \tilde\mg$ 
are then obtained as classically; see in particular (\ref{large_cr}).
From Section 3.1, the classical algebras are
$\tilde B_n = C_n$ and $\tilde C_n = B_n$ if $M$ is even, and 
$\tilde \mg = \mg$ otherwise.
This shows explicitly the refinements of (\ref{frobenius}) which 
arise for even roots of unity, in the most important case of
finite--dimensional representations. Notice that for odd roots of unity,
$K_i$ evaluates to 1 on the special points $\l_z$, 
and Theorem \ref{classical_algebra} essentially reduces to 
(\ref{frobenius}). The general, abstract result is
given in \cite{lusztig_book}.

To summarize the results of this section, 
any $L^{res}(\l)$ for dominant integral $\l$ 
is a direct sum of irreducible \reps of $U_q^{fin}$,
which are related by an action of the classical $U(\tilde \mg)$, extended
by parity generators $K_i$ for even roots of unity. In particular,
this holds for unitary representations. 

\sect{Noncompact forms and unitary representations}

We first recall some concepts in the classical case,
see e.g. \cite{helgason,cornwell}.
Consider a not necessarily compact semisimple Lie group $G$ 
with real Lie algebra $\mg$. 
Let $-\sigma$ be the conjugation on the complexification $\mg_C$ 
with respect to $\mg$ extended as an involution, 
by which we mean an anti-linear anti--algebra 
map whose square
lis the identity; one could equally well consider algebra maps.
On the other hand, the compact form $\mg_K$ of $\mg_C$ is the eigenspace
with eigenvalue $-1$ of the Cartan--Weyl involution $\theta$.
By a theorem of Cartan (see \cite{helgason}, Theorem 7.1),
one can assume that $\sigma = \phi \circ \theta$, 
where $\phi$ is a linear automorphism of $\mg_K$ with $\phi^2=1$.
Let $\mathfrak{k}$ be the eigenspace of $\phi$ with eigenvalue $+1$, and
$\mathfrak{p}$ the eigenspace with eigenvalue $-1$. Then 
$\mg = \mathfrak{k}\oplus i\mathfrak{p}$ is the 
Cartan decomposition of $\mg$, and $\mathfrak{k}$ is a maximal compact
subalgebra. A root $\a$ is called compact if the corresponding 
root vector is in $\mathfrak{k}$.
The star structure is then defined as $\ast = \sigma$.

Now there are two cases, depending on if $\phi$  is an
inner automorphism or an outer automorphism \cite{cornwell}. 
In this work, we only consider the first type, which 
covers $so(n,2p), su(n,p)$, $so^*(2l), sp(n,p)$, $sp(l,\R)$, 
and various forms for the exceptional groups.
We will find quantum versions and unitary \reps
for  all them, even though not all of the \reps will have a 
classical limit.
The second type includes $sl(l+1,\R), su^*(l+1), so(2l-2p-1,2p+1)$,
and exceptional cases.

Up to equivalence, the inner automorphisms 
of a simple Lie algebra of rank $r$
are given by $2^r$  ``chief'' inner automorphisms of the form
$\phi(H_i) = H_i, \; \phi(X_i^{\pm}) = s_i X_i^{\pm}$,
for $s_i = \pm 1$ (\cite{cornwell}, Ch. 14).
They define the real forms
\berr
H_i^*      &=& H_i, \nn\\
(X_i^{\pm})^* &=& s_i X_i^{\mp}, \quad \mbox{for }\; s_i = \pm 1.
\label{s_i_reality}
\err
They are not necessarily inequivalent;
the compact case corresponds to all $s_i=1$. It should however be noted
that real forms which are equilvalent classically are not
necessarily equivalent in the $q$--deformed case. 
For example, the real form $(X^{\pm})^* = -X^{\pm},\; H^* = -H$ for $|q|=1$
of the ''non--restricted'' $U_q(sl(2,\R))$ considered in \cite{schmuedgen} 
is classically equivalent to the form
$(X^{\pm})^* = -X^{\mp}, \; H^* = H$, which is a special case of
(\ref{s_i_reality}). Nevertheless, the first form 
has no unitary \reps at roots of unity if imposed on 
$U_q^{fin}(sl(2))$, while the second does. 

We consider $U_q^{fin}$, which becomes a $*$--algebra for any of the 
forms (\ref{s_i_reality}) for $q$ a root of unity.
Now we allow non--integral weights as well (it should be noted that the 
weights must be integral if working with $U_q^{res}(\mg)$). 
Then there exist one--dimensional \reps
$L^{fin}(\l_r)$ of $U_q^{fin}$ with weight
$\l_{r} = \sum_i r_i M_i \L_{i}$ generalizing (\ref{l_z}), 
where $r_i \in \Q$ such that $[r_i M_i ]_{q_i} =0$, or equivalently
$q^{2r_i M_i d_i} =1$ for all $i$.
This follows immediately from (\ref{keller_eq}). Explicitely, 
\be
\l_r = \sum_i \frac{m}{2nd_i} p_i \L_i
\label{l_r_z}
\ee
with $p_i \in \Z$.

Let $L^{fin}(\l)$ be a unitary \rep of the compact form (such as in Theorem 
\ref{cpct_thm}) with inner product $(\;,\;)$, and consider 
$L^{fin}(\l) \tens L^{fin}(\l_{r})$. This is again an irreducible \rep
of $U_q^{fin}$, and we can define an inner product on it by 
\be
(v\tens \rho_{r},w \tens \rho_{r}):=(v,w)
\label{inner_z}
\ee
where $\rho_{r} \in L^{fin}(\l_{r})$. It is positive
definite by definition. Let us calculate the
adjoint of $X_i^{\pm}$ on this Hilbert space:
\berr
(v\tens \rho_{r},X_i^{\pm}\cdot (w\tens \rho_{r}))
  &=& (v\tens \rho_r, X_i^{\pm}\cdot w \tens q_i^{-r_i M_i/2}\rho_{r}) \nn\\
  &=& q_i^{-r_i M_i/2} (v,X_i^{\pm}\cdot w). 
\err
On the other hand, 
\berr
(X_i^{\mp} \cdot (v\tens \rho_{r}),w\tens \rho_{r})
  &=& (X_i^{\mp}\cdot v \tens q_i^{-r_i M_i/2}\rho_{r},w\tens \rho_{r}) \nn\\
  &=& (q_i^{-r_i M_i/2} v,X_i^{\pm}\cdot w)
\err
by unitarity of $L^{fin}(\l)$. By definition, the inner product is 
antilinear in the first argument. Now there are 2 cases: first, 
if $q_i^{M_i r_i} =1$, then $q_i^{M_i r_i/2} =\pm 1$, and
the adjoint of $X_i^{\pm}$ becomes $(X_i^{\pm})^* = X_i^{\mp}$.
Second, if $q_i^{M_i r_i} = -1$, then $q_i^{M_i r_i/2} =\pm i$, and
the adjoint of  $X_i^{\pm}$ is $(X_i^{\pm})^* = -X_i^{\mp}$.
Therefore we have proved

\begin{theorem}
Let $L^{fin}(\l)$ be a unitary \rep of the compact form of $U_q^{fin}$,
and $L^{fin}(\l_r)$ a one--dimensional \rep of
$U_q^{fin}$ with weight $\l_r = \sum_i r_i M_i \L_i$ as in (\ref{l_r_z}). 
Then $L^{fin}(\l+\l_r) = L^{fin}(\l) \tens L^{fin}(\l_{r})$ 
with inner product (\ref{inner_z})
is a unitary \rep of the real form (\ref{s_i_reality})
of $U_q^{fin}$, where $s_i = q_i^{M_i r_i} = <K_i,\l_{r}> = \pm 1$.
All unitary \reps of that real form can be obtained in this way.
\label{nc_reps_thm}
\end{theorem}

The last statement follows since the noncompact \reps can similarly 
be ``shifted'' back to the compact form.

This explains the role of the extra, ``quasiclassical'' generators $K_i$
at even roots of unity: they determine the real form of a representation.
While the symmetric form of the coproduct (\ref{coproduct_X})
was useful in the proof, it is irrelevant for the result. 

For the remainder of this section we  concentrate on the case of 
integral weights, i.e. $\l_r = \l_z$ as in (\ref{l_z}),
and determine which of the 
classical noncompact forms actually occur in this way.

In the simply laced case, $M_i =M$, and
$q_i^{M_i} = -1$ precisely if $q$ is an even root of unity.
Thus for odd roots of unity, $s_i =1$ for all $i$,  
whereas for even roots of unity, $s_i = (-1)^{z_i}$,
so that there
are unitary \reps for all the noncompact forms considered. 

In the non--simply laced case, consider first $B_n,C_n$ and $F_4$.
If $q$ is odd, i.e. $q_s^{M_s} =1$ with odd $M_s = M$, 
then $M_i = M$, and  $s_i=1$ for all $i$. Therefore 
only the compact form occurs. 
For even $q$, one has to distinguish whether $M=m/2$ is even or odd.
If $M$ is odd, then $M_l = M_s =M$, therefore $q_l^{M_l} =1$ 
and $q_s^{M_s} =-1$.
This means that only those noncompact forms with $s_i=1$ for $\a_i$ a 
long root and $s_i=(-1)^{z_i}$ for $\a_i$ a short root occur.
If $M$ is even, then $M_l = M_s/2$, and $q_i^{M_i} =-1$ for all $i$. 
Therefore $s_i = (-1)^{z_i}$ for all $i$, and again all
noncompact forms considered are realized
(to recover the results in \cite{ads_paper}, notice that the conventions
there are such that $d_s = \frac 12$). 

Finally consider  $G_2$. If $q$ is odd, then $M_l = M/3$ if $M$ is a multiple
of 3, and $M_l = M$ otherwise. In either case, $q_i^{M_i} =1$ for all $i$,
and $s_i=1$ for all $i$.  If 
$q$ is even, then $q_i^{M_i} =-1$ for all $i$, thus 
$s_i = (-1)^{z_i}$ for all $i$, and again all
noncompact forms considered are realized.

The classical limit of these unitary \reps will be discussed  in Section 7.
Notice that Theorem \ref{nc_reps_thm} also yields additional 
unitary \reps of the compact form with generally non--integral weights, 
for $q_i^{M_i r_i}=1$. We will see in Proposition \ref{plane_prop} however
that the distance of their weights from the origin becomes 
infinite as $q$ approaches 1. In that sense, they are non--classical.

\sect{Reality--preserving algebra on $L^{res}(\l)$}

Consider a dominant integral weight $\l_0$ such that
$L^{fin}(\l_0)$ is a unitary \rep of the compact form with
$0 \leq (\l_0,\ad_i) < M_i$ for all $i$, and a special point $\l_z$.
By Theorem \ref{L_res_theorem}, 
$L^{res}(\l_0+\l_z) = \oplus_{z'} (L^{fin}(\l_0) \tens L^{fin}(\l_{z'}))$ 
is a direct sum of irreducible \reps of $U_q^{fin}$, where the
$L^{fin}(\l_{z'})$ are one--dimensional components of $L^{res}(\l_z)$.
These sectors are unitary \reps of various real forms of $U_q^{fin}$, 
according to Theorem \ref{nc_reps_thm}. Moreover by
Theorem \ref{classical_algebra}, the ``large''
generators $\tilde X_i^{\pm}$ connect the various sectors
with different $z'$. It is natural to ask which subalgebra of the
classical $U(\tilde \mg)$ connects only those sectors with the same
real form. This will be called reality--preserving algebra.
Of course, the $(\tilde X_j^{\pm})^2$ always preserve the real form,
but they do not form a closed algebra. 

$X_{\tilde\a}^{\pm}$ as defined below Theorem \ref{classical_algebra} 
preserves the real form if and only if 
\be
[\tilde X_{\a}^{\pm}, K_i] =0 \quad \mbox{for all }\; i,
\ee
or
\be
q^{M_{\a}(\a,\a_i)} =1  \quad \mbox{for all }\; i
\label{inv_alg_cond}
\ee
This is equivalent to $q^{M_{\a}(\a,\b)} =1$ for all roots $\b$.
Using the Weyl group, we can assume that $\a = \a_j$ is a 
simple root, since all other $\a$ satisfying (\ref{inv_alg_cond}) 
are then obtained as the image under the Weyl group of the 
simple ones.

First consider the simply laced case. Then for any $j$, there
is an $i$ such that $(\a_i,\a_j) = -1$, therefore 
$q^{M_j(\a_i,\a_j)} =1$ only if $q$ is odd. But then all sectors are
all compact. Therefore the reality--preserving algebra is 
$\mg$ for odd $q$, and trivial otherwise.

Next consider the non--simply laced case. If $q$ is odd, then
all forms are compact, and the reality--preserving
algebra is clearly $\mg$. 

Thus assume $q$ is even. 
For $G_2$, $q^{M_j(\a_i,\a_j)} = q^{-3 M_j} =-1$ if $i\neq j$,
and the reality--preserving algebra is trivial.

For $B_n,C_n$ and $F_4$, 
observe first that if $A_{ij} \neq 0$ and $d_j \geq d_i$, 
then $(\a_i,\a_j) = -d_j$. Therefore 
$q^{M_j(\a_i,\a_j)} = q^{M_j \max\{d_i,d_j\}}$. 

One has to distinguish $M=m/2$ even and odd.
Assume $M$ is even, so that $M_l = M_s/2 = M/2$. 
Then the only way that $q^{M_j\max\{d_i,d_j\}}=1$ for all $i\neq j$ 
with $A_{ij} \neq 0$ is $M_j = M$ and $\max\{d_i,d_j\} = 2$, 
i.e. $j$ is short and is connected only to  long nodes in the 
Dynkin diagram. The only case where this happens is $B_n$,
which has one short simple root. By the Weyl group, it follows that
the reality--preserving algebra is generated by all $\tilde X_{\a_s}^{\pm}$
where $\a_s$ are the short roots of $B_n$. Since $q$ is even, the dual algebra 
$\tilde \mg$ of $B_n$ is $C_n$, i.e. these $\tilde X_{\a_s}^{\pm}$ 
correspond precisely to 
the long roots of $\tilde \mg$. Now $C_n$ has precisely $n$ 
long roots which are all orthogonal, 
and the corresponding root vectors commute. Therefore for even $M$,
the reality--preserving algebra for $B_n$ is $(su(2))^n$, generated by the
$\tilde X_{\a_s}^{\pm}$ which commute with each other (on $L^{res}(\l)$).
For $C_n$ and $F_4$, it is trivial except for $C_2 \cong B_2$.

If $M$ is odd, then $M_l= M_s = M$, thus $q^{M d_l} =1$,
and $q^{M d_s} =-1$. Therefore if
$q^{M \max\{d_i,d_j\}}=1$ for all $i\neq j$ 
with $A_{ij} \neq 0$,  then either $j$ must be long, or  
$j$ is connected only to long nodes in the Dynkin diagram. 
For $B_n$ this holds for all $j$, for $C_n$ this holds for the 
one long simple root, and for $F_4$ this holds for the 2 long simple roots.
Therefore the reality--preserving algebra for $B_n$ is again $B_n$,
with generators $\tilde X_{\a}^{\pm}$ for all $\a$. For $C_n$,
it is $(su(2))^n$ with generators $\tilde X_{\a_l}^{\pm}$ 
which commute with each other, where $\a_l$ are the long roots of $C_n$.
For $F_4$, it is the algebra generated by all long roots,
which is $D_4$.

\sect{Non--integral weights and the classical limit}

In this section, we want to determine which of the unitary \reps
of $U_q^{fin}$ in Theorem \ref{nc_reps_thm} have a well--defined
classical limit. The idea is to consider them as highest--weight
modules in a suitable way with fixed highest weight, 
and let $q$ approach 1.

For dominant integral $\l_0$ and 
$\l_r = \sum p_j\frac{m}{2nd_{j}} \L_j$ as in (\ref{l_r_z}), consider 
\berr
L^{fin}(\l) &=& L^{fin}(\l_0) \tens L^{fin}(\l_{r}) \qquad {\rm with} 
                    \nonumber\\
(\l_0+\rho,\ad) &\leq& \lceil\frac m{2nd_{\a}}\rceil + 1 \quad
          \mbox{for all }\;  \a \in Q^+,
\label{bounds}
\err
where  $\l=\l_0+\l_r$. 
According to Theorems  \ref{cpct_thm} and  \ref{nc_reps_thm}, this 
is a unitary \rep of a certain noncompact form determined by $\l_r$.
We want to understand the location
of the weights of $L^{fin}(\l)$ in weight space, and in particular
if they are close enough to the origin so that they can have 
a classical limit as a highest weight module. For $q \neq 1$ of course, 
they can always be viewed as highest weight modules.

The bound (\ref{bounds}) for $L^{fin}(\l)$ being unitary
can  be stated more geometrically as follows.
Divide weight space into alcoves separated by the hyperplanes 
\be
h_{\a}^z:=\{\mu; \; (\mu,\a) = \frac m{2n}z \}
\label{h_planes}
\ee 
for all roots $\a$ and $z \in \Z$, 
similar as in Section 3.1. Then $\l_0$ is 
in the fundamental alcove by (\ref{bounds}), using the fact that
$(\rho, \ad) \geq 1$ for all positive roots $\a$;
the latter can be seen using $\rho = \sum_i \L_i$.
By the Weyl group, all weights of $L^{fin}(\l_0)$ are therefore 
contained in the union of those alcoves which have the origin 
as corner, more precisely within a certain distance from its walls as 
determined by (\ref{bounds}). Since the set of
hyperplanes (\ref{h_planes}) is invariant under 
translations by $\l_r$, 
the weights of $L^{fin}(\l)$ are contained in the union
of those alcoves which have $\l_r$ as corner.
In particular, they are contained in a half--space with the 
origin on its boundary. Since the distance between parallel hyperplanes 
goes to infinity as $q\rightarrow 1$, 
$L^{fin}(\l)$ can have a classical limit only if $\l_r$ and the origin 
belong to the same alcove. This puts a restriction on the possible 
real forms as determined by $\l_r$.

To make this more precise, recall the definition of compact roots  
in Section 5, and the definition of the Coxeter labels $a_i$ 
which are the coefficients
of the highest root $\theta = \sum_j a_j \a_j$, and satisfy 
$a_i \geq 1$ for all $i$. 

\begin{prop}
Consider $\l_r = \sum p_j\frac{m}{2nd_{j}} \L_j$ with $p_j \in \Z$, 
as in (\ref{l_r_z}). It belongs to the closure of an alcove as defined above
which also contains the origin as a corner, if and only if there is a set of 
simple roots (denoted again by $\a_i$) such that $\a_{i_0}$
has Coxeter label $a_{i_0}=1$ and $p_j = \pm\d_{j,i_0}$, hence
\be
\l_r = \pm \frac{m}{2nd_{i_0}} \L_{i_0}.
\label{l_r}
\ee
In that case, the remaining $r-1$ simple roots are 
compact with respect to the real form defined in Theorem (\ref{nc_reps_thm}).
\label{plane_prop}
\end{prop}

\begin{proof}
By using an element $\sigma$ of the Weyl group if necessary, we can assume
that $\l_r$ is an anti--dominant weight.
Then $\l_r = -\sum p_j\frac{m}{2nd_{j}} \L_j$ with $p_j \in \N$, and
$(\l_r,\theta) = -\sum p_j a_j\frac m{2n}$.
This implies that $(\l_r,\theta) \leq - \frac m{2n} = (h^{-1}_{\theta},\theta)$
where $h^{-1}_{\theta}$ is a hyperplane as defined in (\ref{h_planes}),
and equality holds precisely if 
$\l_r = -\frac{m}{2nd_{i_0}} \L_{i_0}$ and  $a_{i_0}=1$ for 
some $i_0$. The desired set of 
simple roots is obtained by applying $\sigma$ to the original  
simple roots. Then all (new) simple roots $\a_j$ for $j \neq i_0$
are compact according the definition in Theorem \ref{nc_reps_thm}.
\end{proof}



We will always use this set of simple roots from now on.
The corresponding generators in $U_q^{res}(\mg)$ could be obtained via the 
braid group action \cite{lusztig}, but this is not needed since
we will only make statements about the characters below. 
The corresponding real form is 
\berr
(X_{i_0}^{\pm})^* &=& - X_{i_0}^{\mp}, \qquad \mbox{and} \nn\\
(X_j^{\pm})^*     &=&  X_j^{\mp} \qquad \mbox{for } \;\; j \neq i_0.
\label{herm_pair}
\err
In the classical limit,
the center of $\mathfrak{k}$ is then one--dimensional and generated by
an element of the Cartan subalgebra dual to 
$\L_{i_0}$, which is orthogonal to the compact roots.
Explicitly, this leads to the following cases:  
\begin{itemize}
\item 
$i_0=1,2,...,l$ for $A_l$, corresponding to $su(l+1-p,p)$ for all $p$
\item 
$i_0=1$ for $B_l$, corresponding to $so(2l-1,2)$
\item
$i_0=1$ for $D_l$, corresponding to $so(2l-2,2)$
\item
$i_0=l$ or equivalently $i_0=l-1$ for $D_l$, corresponding to $so^*(2l)$
\item
$i_0=l$ for $C_l$, corresponding to $sp(l,\R)$
\item
$i_0=1$ or equivalently $i_0=5$ for $E_6$, and $i=6$ for $E_7$,
\end{itemize}
see for example \cite{cornwell}, table 14.1.
Not surprisingly, these are precisely the cases
where highest weight modules exist in the classical limit, see 
\cite{enright,garland} and references therein.
We will restrict ourselves to (\ref{herm_pair}) from now on, and show 
how to recover the classical unitary highest weight \reps 
from the $L^{fin}(\l)$. To do that,
we choose a minus sign in (\ref{l_r}), since then
it is possible to consider
highest--weight modules $L^{fin}(\l)$ with fixed highest weight $\l$
independent of $q$, in  particular as $q \rightarrow 1$. 
The plus sign would correspond to lowest--weight 
modules in the classical limit. 

To make the connection with the literature 
on the classical case \cite{enright}, 
consider the character 
$\chi(L(\l+z\L_{i_0}))e^{-z\L_{i_0}}$ for $z\in \R$, where 
$L(\l+z\L_{i_0})$ is the classical irreducible highest weight 
module with highest weight $\l+z\L_{i_0}$.
It is independent of $z$ for sufficiently negative $z$, which can be seen
from the strong linkage principle (see e.g. \cite{kac}): by writing
$\l = c_0 \L_{i_0} + \sum_{j \neq i_0} n_j \L_j$ and noticing that
the compact roots are orthogonal to $\L_{i_0}$, if follows that
for sufficiently negative $z$, all
weights strongly linked to $\l+z\L_{i_0}$ are in 
the orbit of the compact Weyl group acting on $\l+z\L_{i_0}$.
The {\em first reduction point} $z_0$ is 
the maximal value of $z$ where this is no longer the case.
Clearly $L(\l+z\L_{i_0})$ can only be unitary with respect to  
(\ref{herm_pair}) if $\l$ is a dominant integral 
weight with respect to $\mathfrak{k}$, i.e. $n_j \in \N$ in
the above notation. Provided this is the case, 
$L(\l+z\L_{i_0})$ is unitary \cite{enright} if and only if
$z\leq z_0$, or $z$ is in a certain finite set of $z>z_0$.

In the $q$--deformed case, we can show the following:
\begin{prop}
Let $\l$ be a rational weight which is dominant integral with 
respect to $\mathfrak{k}$. If the first reduction point of 
$L(\l+z\L_{i_0})$ is at $z\geq 0$ for $q=1$, then there exists 
a series of roots of unity  $q_k \rightarrow 1$ such that
$L^{fin}(\l)$ is unitary with respect to 
(\ref{herm_pair}) for all $q=q_k$. In particular, this holds if 
$(\l+\rho,\ad) \leq \lceil(-\l,\ad_{i_0})\rceil + (\l,\ad_{i_0}) +1$
for all positive noncompact roots $\a$.
\end{prop}
Of course, this generalizes to irrational $\l$ which can be approximated by 
rational weights as above.
\begin{proof}
Assume that $\l$ is as required. Then  
there are $m,n\in \N$ such that $(\l-\l_r,\ad_{i_0}) \in \N$
where $\l_r = -\frac{m}{2nd_{i_0}} \L_{i_0}$, thus 
$\l-\l_r$ is dominant integral. For $k\in \N$, let
$m_k:=m+2nkd_{i_0}$, $\l_{r,k}:=-\frac{m_k}{2nd_{i_0}} \L_{i_0}$,
and $q_k:=e^{2\pi i n/m_k}$. We claim that for sufficiently large $k$,
the character of $L^{fin}(\l-\l_{r,k})$ 
is given by Weyls formula for all $q'$ between 1 and $q_k$. 
Then the first part of the
proposition follows from
Theorems \ref{cpct_thm} and \ref{nc_reps_thm}.

Let $q'=e^{2\pi i n'/m'}$ with $\frac{n'}{m'} < \frac nm$,
with associated $M_{\a}'$ as in Section 2.
By the strong linkage principle \cite{anderson,ads_paper}, the
character of $L^{fin}(\l-\l_{r,k})$ can differ from $\chi(\l-\l_{r,k})$
only by the sum of classical characters $\chi(\mu_n)$ (\ref{char})
with dominant $\mu_n$, which are ``strongly linked'' 
to $\l-\l_{r,k}$ by a series of reflections by hyperplanes 
${H_{\a}^z}'$ defined as in (\ref{H_z}) using $M_{\a}'$, but 
shifted by $-\rho$. They again divide weight space into (shifted) alcoves,
with corresponding special points for $q'$, also shifted by $-\rho$. 

Now $a_{i_0}=1$ implies that $-\rho$ and $-\l_{r,k} - \rho$ 
are in the same shifted alcove for any $q'$ between 1 and $q_k$,
because $(-\l_{r,k},\ad) \leq \frac{m_k}{2nd_{\a}} \leq M_{\a}'$
for all positive $\a$. 
Moreover, the union of the alcoves which have $-\rho$ as a 
corner is a convex set of weights, and invariant under the Weyl group
action with center $-\rho$. Therefore all weights in that set which
are strongly linked to $-\l_{r,k}-\rho$ are obtained by the action of the 
{\em classical} Weyl group with center $-\rho$; in particular, $-\rho$ 
is not. Thus if $k$ is large enough, all {\em dominant} $\mu_n$ 
strongly linked to $\l-\l_{r,k}$ can be obtained by reflections 
of $\l-\l_{r,k}$ by those
${H_{\a}^z}'$ which contain the special point $M_{i_0}' \L_{i_0}-\rho$.
However using the assumption, the character of $L^{fin}(\l-\l_{r,k})$ 
is not affected by these $\mu_n$: 
indeed, by a shift by $-\frac{m'}{2n'd_{i_0}}\L_{i_0}$ as in Section 5, 
$L^{fin}(\l-\l_{r,k})$ can be related to $L^{fin}(\l+z\L_{i_0})$
for $z = \frac{m}{2nd_{i_0}}-\frac{m'}{2n'd_{i_0}} <0$. 
The special point $M_{i_0}' \L_{i_0}-\rho$ is then moved to 
$(M_{i_0}'-\frac{m'}{2n'd_{i_0}})\L_{i_0}-\rho$, and is only relevant 
for large $k$ if $M_{i_0}' = \frac{m'}{2n'd_{i_0}}$, 
when it becomes $-\rho$. However by the
assumption on the first reduction point, the character of the 
classical $L(\l+z\L_{i_0})$ is not affected by 
the hyperplanes through $-\rho$
for $z<0$. Using the fact that $U_q(\mg)$ is the same as $U(\mg)$
as algebra over $\C[[q-1]]$  \cite{drinfeld_2}, 
the character of $L^{fin}(\l+z\L_{i_0})$ is not affected by 
the hyperplanes through $-\rho$ either. Combining all this,
it follows that the character of $L^{fin}(\l-\l_{r,k})$ is given 
by Weyls formula for all $q'$ between 1 and $q_k$.

In particular, this holds if 
$(\l-\l_{r,k}+\rho,\ad) \leq \lceil\frac{m_k}{2nd_{\a}}\rceil+1$
for all positive noncompact roots $\a$, by Theorem \ref{cpct_thm}.
This is certainly satisfied for compact $\a$ if $k$ is sufficiently large.
Using $\frac{m_k}{2nd_{i_0}} = -(\l,\ad_{i_0}) +(\l-\l_{r,k},\ad_{i_0}) 
\in -(\l,\ad_{i_0}) +\N$ and the fact that $d_{i_0} =2$ for the 
non--simply laced cases, this bound follows from the given condition.
\end{proof}

Therefore we recover the classical results on unitary 
highest weight representations, except for the small, 
finite set of $z> z_0$, which we cannot address here. 
It is quite possible that there exist  
unitary representations of $U_q^{fin}$ corresponding to these remaining 
cases; this would have to be studied by other methods. 
By Theorem \ref{nc_reps_thm}, they would correspond to additional 
unitary \reps of the compact form, as was pointed out in Section 3.

In the example of the Anti--de
Sitter group, the highest--weight property 
corresponds to positivity of the energy
\cite{ads_paper}, which is an important physical requirement.

Notice that the unitary \reps of noncompact forms of $U_q^{fin}$ 
in general have non--integral, but rational weights.
Those with integral weights can have a classical limit
only if $\frac m{2nd_{i_0}} \in \Z$, in particular
$q$ must be an even root of unity.

The question arises if and how the unitary \reps 
of $U_q^{fin}$ in those cases where there exists no
classical unitary highest weight \rep might be related to other 
classical series of unitary
representations, and how the latter may be obtained from
the quantum case at roots of unity. The answer may be related to the
fact that there do exist other types of unitary \reps 
of the non--restricted specialization for $|q|=1$,
such as $U_q(sl(2,\R))$ \cite{schmuedgen}, 
as was mentioned in Section 5. This certainly deserves
further investigation.

\sect{Acknowledgements}
I would like to thank Konrad Schm\"udgen and David Vogan
for pointing out some related references, and John Madore for 
useful comments.

\begin{appendix}
\section{Appendix}

We prove Theorem \ref{classical_algebra} by 
verifying the classical relations of the Chevalley basis 
$\tilde X_i^{\pm}$ and $\tilde H_i$.

To calculate $\tilde H_i$ on the weights $\l_{z} = \sum_j z_j M_j\L_j$, 
one can use the standard commutation relation 
\be
[X_i^{+(M_i)}, X_i^{-(M_i)}] = 
    \left[\begin{array}{c} H_i\\ M_i \end{array}\right]_{q_i} + u_q^{fin},
\label{id_1}
\ee
where the last term vanishes on $L^{res}(\l_{z})$. We also need
the following identity \cite{lusztig_book} which can be checked directly: 
If $q^{2M} =1$, then
\be
\left[\begin{array}{c} aM+b\\ cM+d \end{array}\right]_q = 
                              q^{M^2c(a+1)+M(ad-bc)}
    \left[\begin{array}{c} b\\ d \end{array}\right]_q 
    \(\begin{array}{c} a\\ c \end{array}\)_1 .
\label{binom_id}
\ee
Furthermore $[K_i,X_i^{\pm (M_i)}] =0$, and therefore
$[\tilde X_i^{+}, \tilde X_i^{-}] = [X_i^{+(M_i)}, X_i^{-(M_i)}] \tilde K_i
              q_i^{M_i^2}$.
Using (\ref{id_1}), (\ref{binom_id}) and (\ref{K_i_eval}), 
this evaluates on $v_{\l_{z'}}$ to
\be
\tilde H_i\cdot v_{\l_{z'}} = 
    \left[\begin{array}{c} z'_i M_i \\ M_i \end{array}\right]_{q_i} 
     q_i^{z'_i M_i^2} q_i^{M_i^2}\cdot v_{\l_{z'}} 
   = z'_i v_{\l_{z'}},
\ee
as claimed.

We next check that
\be
[\tilde X_i^{+}, \tilde X_j^{-}] = 0
\ee
for $i\neq j$. This is clear if $A_{ij} =0$. Otherwise, one can write
\be
X_i^{+(M_i)} \tilde K_j = s_{ij} \tilde K_j X_i^{+(M_i)},
\ee
where 
$s_{ij} = q_i^{M_i M_j A_{ji}} = q^{M_i M_j (\a_i,\a_j)} = s_{ji} = \pm 1$. 
Then $s_{ji}^{a_i} s_{ij}^{1-a_j} = q^{M_i M_j(\a_i,\a_j)2a_i} =1$,
since $a_i = (1-a_j)$ if $A_{ij} \neq 0$ and $i\neq j$. Therefore
\be
\tilde X_i^{+} \tilde X_j^{-} = s_{ji}^{a_i} s_{ij}^{1-a_j}
\tilde X_j^{-} \tilde X_i^{+} = \tilde X_j^{-} \tilde X_i^{+}.
\ee

Next, to verify  
\be
[\tilde H_i, \tilde X_j^{\pm}] = \pm \tilde A_{ji} \tilde X_j^{\pm},
\ee
replace again $\tilde H_i$ by 
$\left[\begin{array}{c} H_i\\ M_i \end{array}\right]_{q_i} \tilde K_i
q_i^{M_i^2}$, and observe using (\ref{dual_cartan}) that
\be
[H_i]_{q_i}X_j^{\pm(M_j)} = X_j^{\pm(M_j)} [H_i\pm M_j A_{ji}]_{q_i} = 
 X_j^{\pm(M_j)} [H_i\pm M_i \tilde A_{ji}]_{q_i}.
\ee 
We first show  
\be
\(\left[\begin{array}{c} H_i\\ M_i \end{array}\right]_{q_i} \tilde K_i
q_i^{M_i^2}\) X_j^{+(M_j)} \tilde K_j^{a_j} = 
X_j^{+(M_j)} \tilde K_j^{a_j}
\(\left[\begin{array}{c} H_i\\ M_i \end{array}\right]_{q_i} \tilde K_i
q_i^{M_i^2} + \tilde A_{ji}\).
\ee
Using the above, this becomes
\be
X_j^{+(M_j)} \tilde K_j^{a_j}
  \left[\begin{array}{c} H_i+\tilde A_{ji} M_i\\ M_i\end{array}\right]_{q_i} 
  \tilde K_i  s_{ji} q_i^{M_i^2} = 
X_j^{+(M_j)} \tilde K_j^{a_j} 
\(\left[\begin{array}{c} H_i\\ M_i \end{array}\right]_{q_i} \tilde K_i
q_i^{M_i^2} + \tilde A_{ji}\).
\ee
Restricting on a weight $\l_{z'}$, it remains to show
\be
\left[\begin{array}{c} (z'_i+\tilde A_{ji}) M_i\\ M_i\end{array}\right]_{q_i} 
  q_i^{M_i^2 z'_i} s_{ji} q_i^{M_i^2} = 
\(\left[\begin{array}{c} z'_i M_i\\ M_i \end{array}\right]_{q_i} 
q_i^{M_i^2 z'_i} q_i^{M_i^2} + \tilde A_{ji}\).
\ee
Now $s_{ji} = q_i^{M_i^2 \tilde A_{ji}}$, and the claim follows from
(\ref{binom_id}). The calculation for 
$[\tilde H_i, \tilde X_j^{-}] = - \tilde A_{ji} \tilde X_j^{-}$
is completely analogous.

Finally, the Serre relations are
\be
[\tilde X_i^{+}, ... ,[\tilde X_i^{+},\tilde X_j^{+}]...]^
{1-\tilde A_{ji}} =0
\ee
($1-\tilde A_{ji}$ brackets) on $L^{res}(\l_{z})$,
and similarly for the negative roots. 

To prove this, consider  a P.B.W. basis of $U_q^{+res}(\mg)$, 
which is given by the expressions
$X_{\a_{N}}^{+(t_{N})}...  X_{\a_{1}}^{+(t_{1})}$ where 
$\{\a_1, ... \a_{N} \}$ is an ordered basis of the positive
roots, obtained e.g. by the braid group action \cite{lusztig_book}.
Let $\tilde Q = \{\tilde \a =  M_{\a} \a\}$ be the set of roots of the 
lattice of special points. For $k\in\N$ such that 
$\tilde \b = k\tilde \a_i+\tilde\a_j \in \tilde Q$, define
$\tilde X_{\tilde\b}^+:=X_{\b}^{+(M_{\b})} \tilde K_i^{k a_i} 
\tilde K_j^{a_j}$
generalizing (\ref{class_generators}),
and $\tilde X_{\tilde\b}^+:=0$ if $\tilde \b\notin \tilde Q$.
We claim that  
\be
[\tilde X_i^+,\tilde X_{\tilde\b}^+] = 
c \tilde X_{\tilde\a_i+\tilde\b}^+
\label{large_cr}
\ee
if acting on $L^{res}(\l_{z})$, 
for some constant $c$. This clearly implies the Serre relations. 
The proof is by induction on $k$, using the well--known commutation relations
\cite{CH_P}
\be
X_{\a_r}^+ X_{\a_s}^+ - q^{(\a_r,\a_s)} X_{\a_s}^+ X_{\a_r}^+ = 
    \sum c(t_{r+1},...,t_{s-1}) 
X_{\a_{s-1}}^{+t_{s-1}}...  X_{\a_{r+1}}^{+t_{s-1}}
\label{XX_cr}
\ee
for $r<s$, with some constant $c(t_{r+1},...,t_{s-1})$.

We want to order the expression $\tilde X_{\tilde\b}^+\tilde X_i^+$
(or the reversed form)
as in the P.B.W. basis, using (\ref{XX_cr}). The leading term is 
\be
q^{M_i M_{\b}(\a_i,\b)}s_{ij}^{a_j}s_{ji}^{a_i}
                           \tilde X_i^+\tilde X_{\tilde\b}^+, 
\label{stuff}
\ee 
since $[\tilde K_i^{a_i},X_i^{+(M_i)}]=0$.
We claim that the only other term on the rhs of (\ref{XX_cr}) 
which may not vanish on $L^{res}(\l_{z})$ is proportional to
$\tilde X_{\tilde\a_i+\tilde\b}^+$. This is so
because only products of ``large'' generators $X_{\a}^{+(M_{\a})}$
are nonzero on $L^{res}(\l_{z})$, and
in fact only one ``large''  generator can occur on the 
rhs of (\ref{large_cr}), because only a simple (formal) pole in $q$
can arise by the derivation property mentioned in Section 2.
Moreover using $\tilde \b = M_{\b}\b = k M_i \a_i + M_j\a_j$, 
it follows that
$q^{M_i M_{\b}(\a_i,\b)} = q^{2kM_i^2 d_i} q^{M_i M_j(\a_j,\a_i)} = 
q^{M_i M_j(\a_j,\a_i)} = s_{ij}$. Thus the overall coefficient
in front of (\ref{stuff}) is $s_{ij}^{1-a_j} s_{ji}^{a_i}$, which is 1 as
above. This concludes the proof.
\end{appendix}


\begin{thebibliography}{9}
 
\bibitem{maldacena} J. Maldacena, ``The Large N Limit of Superconformal 
  Field Theories and Supergravity'' {\em Adv.Theor.Math.Phys.} {\bf 2}, 
   231 (1998) 
\bibitem{drinfeld} V. Drinfeld, "Quantum Groups" {\em Proceedings
  of the International Congress of
  Mathematicians, Berkeley, 1986} A.M. Gleason (ed.), p. 798, AMS,
      Providence, RI
\bibitem{FRT}L.D.Faddeev, N.Yu.Reshetikhin, L.A.Takhtajan.
  "Quantization of Lie
  Groups and Lie Algebras"  {\em Algebra Anal.} {\bf 1} 178 (1989)
\bibitem{jimbo} M. Jimbo, "A q -- Difference Analogue of $U(\mg)$ and the
   Yang -- Baxter
 Equation" {\em Lett. Math. Phys} {\bf 10}, 63 (1985)
\bibitem{klimyk} V. A. Groza, N. Z. Iorgov, A. U. Klimyk, 
  ``Representations of quantum algebra $U_q(u_{n,1})$'' math/9805032;
   Klimyk A. and S. Pakuliak, ``Representations of the quantum algebras
   $U_q(u_{r,s})$ and $U_q(u_{r+s})$ related to the quantum hyperboloid and
   sphere'' {\em J.Math. Phys.} vol. {\bf 33}, No. 6, 1987 (1992)
\bibitem{guizzi} V. Guizzi, ``A classification of unitary highest 
  weight modules of
  the quantum analogue of the symmetric pairs $(A_n,A_{n-1})$'', {\em J. Alg.}
  vol. {\bf 192}, 202 (1997)
\bibitem{ads_paper} H. Steinacker, "Finite--dimensional Unitary 
  Representations of quantum Anti--de Sitter Groups at Roots of Unity"
  {\em Comm. Math. Phys.} {\bf 192}, 687 (1998)
\bibitem{dobrev} V. K. Dobrev, P. J. Moyan, {\em Phys. Lett.} 
   {\bf 315B}, 292 (1993);
     L. Dabrowski, V.K. Dobrev, R. Floreanini and V. Husain, ``Positive 
     energy representations of the conformal quantum algebra'', 
     {\em Phys.Lett.} {\bf B302}, 215-222 (1993)
\bibitem{flato} M. Flato, L.K. Hadjiivanov, I.T. Todorov,
 "Quantum Deformations of Singletons and of Free Zero-Mass Fields"
  {\em Foundations of Physics}, vol.{\bf 23} (4), 571-586 (1993)
\bibitem{schmuedgen} K. Schm\"udgen, ``Operator \reps of 
  ${\cal U}_q(sl_2(\R))$'',  {\em Lett. Math. Phys}. {\bf 37}, 211 (1996)
\bibitem{lusztig_90} G. Lusztig,  "Quantum Groups at roots of 1"
   {\em Geom. Ded.} {\bf 35}, 89 (1990)
\bibitem{lusztig_book}  G. Lusztig,  "Introduction to Quantum
  Groups". {\em Progress in Mathematics} Vol. 110, Birkhaeuser 1993
\bibitem{lusztig} G. Lusztig, "On quantum groups" 
   {\em J. Algebra} {\bf 131}, 466 (1990)
\bibitem{lusztig_mod} G. Lusztig, "Quantum deformations of certain simple
  modules over enveloping algebras", {\em Adv. in Math.} {\bf 70}, 237 (1988)
\bibitem{rosso} M. Rosso, "Finite Dimensional Representations of the Quantum 
  Analog of the Enveloping Algebra of a Complex Simple Lie Algebra" 
   {\em Comm. Math. Phys.} {\bf 117}, 581 (1988)
\bibitem{CH_P} V. Chari and A. Pressley, "A guide to quantum
  groups".   Cambridge University press, 1994
\bibitem{anderson} H.H. Anderson, P. Polo, W. Kexin, {\em Invent. math.}
  {\bf 104}, 1 (1991)
\bibitem{helgason} S. Helgason, ``Differential Geometry,
  Lie Groups and Symmetric Spaces''. Academic Press, New York 1978
\bibitem{cornwell} J.F. Cornwell, ``Group Theory in Physics'', Vol.II.
  Academic Press 1984  
\bibitem{enright} T. Enright, R. Howe, N.R. Wallach, ``A classification 
 of unitary highest weight modules'', in {\em Representation Theory of 
 Reductive Groups}, T. Trombi (Ed.), Progress in Mathematics,
 Birkh\"auser, Boston 1982.
\bibitem{garland} H. Garland, G.J. Zuckerman, ``On unitarizable highest 
  weight modules of Hermitian pairs''. {\em J. Fac. Sci. Univ. Tokyo Sect. 
  IA Math.} {\bf 28}, no. 3, 877 (1982)
\bibitem{kac} V. Kac, D. Kazhdan, "Structure of Representations with 
  Highest Weight
  of infinite-dimensional Lie Algebras", {\em Adv. Math.} {\bf 34}, 97 (1979)
\bibitem{drinfeld_2} V. Drinfeld, "On Almost Cocommutative Hopf Algebras"
   {\em Leningrad Math. J.} {\bf 1}, No. 2, 321 (1990) 


\end{thebibliography}
\end{document}